\documentclass[11pt]{article}
\usepackage{amssymb,amsfonts,amsmath,mathrsfs}

\newtheorem{theo}{Theorem}[section]

\newtheorem{coro}[theo]{Corollary}

\numberwithin{equation}{section}

\def\pf{\noindent {\it Proof.} }
\def\qed{\hfill \rule{4pt}{7pt}}
\begin{document}
\title{The Algorithm Z and  Ramanujan's $_1\psi_1$ Summation}
\author{Sandy  H.L. Chen$^1$, William Y.C. Chen$^2$, Amy M. Fu$^3$ and Wenston J.T. Zang$^4$\\ \ \\
Center for Combinatorics, LPMC-TJKLC\\
Nankai University, Tianjin 300071, P.R. China\\ \ \\
 $^1$chenhuanlin@mail.nankai.edu.cn, $^2$chen@nankai.edu.cn,
$^3$fu@nankai.edu.cn,\\
 $^4$wenston@cfc.nankai.edu.cn}
\date{}
\maketitle
\begin{abstract}
We use the Algortihm Z  on partitions due to Zeilberger, in a variant form,
 to give a  combinatorial proof of Ramanujan's
$_1\psi_1$ summation formula.
\end{abstract}



\section{Introduction} \label{S:P*}

Ramanujan's sum for $_1\psi_1$ has been extensively studied in the
theory of $q$-series,  which is usually stated in the following
form:
\begin{equation}
 _1\psi_1 (a; b; q,
z)=\sum_{n=-\infty}^{\infty}\frac{(a; q)_n}{(b; q)_n}z^n=\frac{(q,
b/a, az, q/az; q)_{\infty}}{(b, q/a, z, b/az;
q)_{\infty}},\,\,|b/a|<|z|<1, \,\,|q| < 1,
\end{equation}
where  the $q$-shifted factorial is defined by
$$
(a;q)_{\infty}=\prod_{n=0}^\infty (1-aq^{n}), \,\, (a;q)_n=(a;q)_{\infty}/(aq^n;q)_{\infty}.
$$

The main result of this paper is a combinatorial proof of the above
formula by using a variation of the Algorithm Z named after
Zeilberger \cite{BressoudZeilberger}. Since Hahn and Jackson
published the first proofs in 1949 and 1950, many other proofs have
been found, see, for example,   Andrews \cite{Andrews1969}, Andrews
and Askey \cite{AndrewsAskey}, Berndt \cite{Berndt},  Fine \cite{Fine}, Ismail
\cite{Ismail}, Mimachi \cite{Mimachi}. However, the combinatorial
proofs have appeared  only recently. Using the Frobenius notation
for overpartitions, Corteel and Lovejoy \cite{CorteelLovejoy} have
found a bijective proof of the constant term identity for the the
following formulation of Ramanujan's $_1\psi_1$ summation:
\begin{equation} \label{yi}
\frac{(-aq;q)_\infty(-bq;q)_\infty}{(q;q)_\infty(abq;q)_\infty}\sum
_{n=-\infty}^\infty\frac{(-a^{-1};q)_n(zqa)^n}{(-bq;q)_n}=\frac
{(-zq;q)_\infty(-z^{-1};q)_\infty}{(bz^{-1};q)_\infty(azq;q)_\infty}.
\end{equation}
 Corteel \cite{Corteel} went on to find a bijection, by using particle
seas, to show  that the coefficients of $z^N$ ($N\neq 0$) on both
sides of \eqref{yi} are equal as well, which leads to the completion
of the combinatorial proof of (\ref{yi}). In the meantime, Yee
\cite{Yee} also found  a  combinatorial proof of \eqref{yi} in the
language  of $F$-partitions defined as pairs of overpartitions with
different sizes written in the two line notation.

In this paper, we shall present a new combinatorial proof of
Ramanujan's $_1\psi_1$ sum based on  a variation of the Algorithm Z.
Conceptually, our bijection is rather simple despite that there are several
steps which do not seem to be avoidable to accomplish the task of
transformations of partitions. As will be seen, the Algorithm Z serves as the
main ingredient of our combinatorial construction for
 Ramanujan's formula. To be precise, our
bijection  is devised for following restatement of Ramanujan's formula
\begin{equation}\label{32}
 \frac{(-q/a; q)_{\infty}(-b/az; q)_{\infty}}{(q;
q)_{\infty}}\sum_{n=-\infty}^{\infty}\frac{(-a; q)_n}{(b; q)_n}z^n=
\frac{(-b/a; q)_{\infty}(-az; q)_{\infty}(-q/az; q)_{\infty}}{(b;
q)_{\infty}(z; q)_{\infty}}.
\end{equation}

The Algorithm Z, as called by Andrews and Bressoud \cite{AndrewsBressoud},
was  found by Zeilberger \cite{BressoudZeilberger}
(Proposition 3.1) as a combinatorial interpretation of the
Gauss coefficient ${n \brack k}$ as defined by the following relation
\[
{1\over (q;q)_{i+j}}{{i+j\brack i}} = {1\over (q;q)_{i}
(q;q)_j}.\]
 Using this algorithm, Andrews and Bressoud have found combinatorial
proofs of some classical $q$-identities. The Algorithm Z has also
been employed by Bessenrodt \cite{BESSENRODT} to give a bijective
proof of a theorem of Alladi and Gordon, and  to
give a combinatorial interpretation of  the Lebesgue identity by Fu \cite{fu}.

\section{The Algorithm Z} \label{S:P*}

In this section, we shall give an overview of the Algorithm Z and use it to
give a combinatorial interpretation of $q$-binomial theorem, which
is an important step of our combinatorial proof of Ramanujan's summation
\eqref{32}:
\begin{equation}\label{binomial}
 \sum_{n
\geq 0}\frac{P_n(b, -a)}{(q; q)_n}z^n=\sum_{n \geq
0}\frac{(-a/b;q)_n}{(q; q)_n}(bz)^n=\frac{(-az; q)_{\infty}}{(bz;
q)_{\infty}}.
\end{equation}where the polynomials\[
P_n(b, -a) =\sum_{k=0}^n {{n \brack k}} a^k
q^{\binom{k}{2}}b^{n-k}=
 \left\{
\begin{array}{ll}
(b+a)(b+aq)\cdots (b+aq^{n-1}),  &\quad \mbox{if} \hskip 0.2cm n \geq 1; \\[6pt]
1,
 &\quad \mbox{if} \hskip 0.2cm n =0,
\end{array}\right.
\]are  the Cauchy polynomials  as called in \cite{chen}.

 A partition $\lambda$ of a nonnegative integer with $r$ parts is denoted by
$\lambda=(\lambda_1, \lambda_2, \ldots, \lambda_r)$, where
$\lambda_1\geq  \lambda_2\geq \cdots \geq  \lambda_r\geq 0$. The
number of parts, called the length of $\lambda$, is denoted by
$l(\lambda)$, and the sum of parts, called the weight of $\lambda$,
is denoted by $|\lambda|$. The conjugate of $\lambda$ is denoted by
$\lambda'=(\lambda'_1,\ldots)$, where $\lambda'_i$ is the number of
positive parts of $\lambda$ that are greater than or equal to $i$. The following bijection is call the Algorithm Z.

\begin{theo}  \label{Zeil}
There is a bijection between the set of pairs of partitions
$(\alpha, \beta)$ where $\alpha$ has $s-r$ parts and $\beta$ has $r$
parts, and the set of pairs of partitions $(\mu, \nu)$, where $\mu$
 has $s$ parts and $\nu$  has $r$ parts with each part not exceeding $s-r$.
 We call $\mu$ the insertion
partition and call $\nu$ the record partition.
\end{theo}

\pf Given a partition $\alpha$ with $s-r$ parts, denoted by
$(\alpha_1, \ldots, \alpha_{s-r})$, and a partition $\beta$ with $r$
parts, denoted by $(\beta_1, \ldots, \beta_r)$, we may insert
$\beta$ into $\alpha$ to create a pair of partitions $\mu$ and
$\nu$. The insertion algorithm can be described as the following
recursive procedure.

\begin{itemize}
\item[$\bullet$ ]If $\beta_1 \leq \alpha_{s-r}$, we insert $\beta_1$ into $\alpha$ so that we get a
new partition $(\alpha_1, \alpha_2, \ldots, \alpha_{s-r+1})$, where
$\alpha_{s-r+1}=\beta_1$. Moreover, we use a zero part as a record
of the insertion position.

\item[$\bullet$ ]
If $\beta_1 > \alpha_{s-r}$, we recursively insert $\beta_1 -1$ into
the partition $(\alpha_1, \alpha_2, \ldots, \alpha_{s-r-1})$.
Suppose that the recursive procedure ends up with $\beta_1 -\nu_1$
being inserted, we use a part $\nu_1$ to record the position of
$\beta_1-\nu_1$. Obviously, we have $0 \leq \nu_1 \leq s-r$.
\end{itemize}

Conversely, given a partition $(\alpha_1, \ldots, \alpha_{s-r+1})$
and a number $\nu_1$ with $0 \leq \nu_1 \leq s-r$, we may extract the
part $\beta_1$ from the given partition. It is easy to see that
above procedure is reversible.

After the part $\beta_1$ has been inserted to $\alpha$, we may
iterate the above procedure to insert remaining parts of $\beta$.
Eventually, we obtain a pair of partitions $(\mu,\nu)$. This
completes the proof. \qed

As an example, taking $\alpha=(5,3,2,1)$, $\beta=(4,3,0)$ with $s=7,
r=3$, we have
$$\mu=( 5, 3, 2,2,2,1, 0),\quad  \nu=(2, 1, 0).$$
Below is the illustration of the insertion procedure
\begin{center}
\begin{tabular}{ccccccccc}
5 & 3 &  & 2 &  &1 & \\
\hline 5 & 3 & 2 & 2 & 2 & 1& 0\\
  & & 2& & 1& &0  \\
  \hline
 & & 4& & 3&  &0\\
\end{tabular}
\end{center}

\begin{coro}\label{4} There is a bijection $\phi$ between the set of
pairs of partitions $(\alpha, \beta)$  and the set of pairs of
partitions $(\mu, \nu)$  satisfying the following conditions
\begin{itemize}
\item[${\bullet}$] $\alpha$ has $i$ distinct parts, $\beta$ has $j$
parts;
\item[${\bullet}$] $\mu$ has $i+j$ parts, $\nu$ has $i$ distinct parts with each part $\leq i+j-1$;
\item[${\bullet}$] $|\alpha|+|\beta|=|\mu|+|\nu|$.
\end{itemize}
\end{coro}

\pf Given a pair of partitions $(\alpha, \beta)$, where $\alpha$ has
$i$ distinct parts and $\beta$ has $j$ parts, we denote by
$\overline{\alpha}$ the partition
$(\alpha_1-i+1,\alpha_2-i+2,\ldots,\alpha_i-0)$. Applying the Algorithm
Z to ($\beta$, $\bar{\alpha}$) yields the desired partition $\mu$
into exactly $i+j$ parts and a partition $\bar{\nu}$ into exactly
$i$ parts with each part  $\leq j$. Set
$\nu=(\bar{\nu}_1+i-1,\ldots,\bar{\nu}_i+0)$. It is clear that
$|\mu|+|\nu|=|\alpha|+|\beta|$. Hence the pair of partitions
$(\mu,\nu)$ satisfy the conditions in the corollary. Since each step
is reversible, we have established a bijection. This completes the
proof.\qed

It is clear that Corollary \ref{4} leads to a combinatorial
proof of the $q$-binomial theorem.  The first partition-theoretic
proof of \eqref{binomial} is due  to Andrews
 \cite{Andrews1967}. There are other proofs of
this classical identity, for example, by overpartitions
\cite{CorteelLovejoy2004} and by  MacMahon diagrams
\cite{JoichiStanton,Pak}.

\section{A Variation of the Algorithm Z}

In this section, we  give  a variation of the Algorithm Z. This algorithm plays a
key role in our combinatorial proof of  Ramanujan's summation
formula.

\begin{theo}\label{7}
 Let $s, t,k,m$ be nonnegative integers. There is
a bijection $\varphi$ between the set of pairs of partitions
$(\alpha, \beta)$ and the set of pairs of partitions $(\mu, \nu)$
satisfying the conditions
\begin{itemize}
\item[$\bullet$ ] $\alpha$ has $s$ distinct parts with each part $\ge m$ and $\beta$ has $t$ parts with each part $\geq k+s+t-1$;
\item[$\bullet$ ] If $s,t>0$, then $\mu$ has $s+t$ distinct nonnegative parts with
$\mu_{s}-\mu_{s+1}\ge m+1$ and $\nu$ has $t$ distinct parts with $ k
\leq \nu_i\leq k+s+t-1$ for each $1\leq i \leq t$;

 If $s>0$ and $t=0$, then $\mu=\alpha$ and $\nu$ is an empty
partition;

 If $s=0$, $t>0$, then
$\mu=(\beta_1-k,\beta_2-k-1,\ldots,\beta_t-k-t+1)$ and
$\nu=(k+t-1,k+t-2,\ldots,k)$;

If $s=t=0$, then both $\mu$ and $\nu$
are empty partitions.
\item[$\bullet$ ] $|\alpha| +|\beta|=|\mu|+|\nu |$.
\end{itemize}
\end{theo}
 \pf Given two partitions $\alpha=(\alpha_1, \ldots, \alpha_{s})$ and
  $\beta=(\beta_1, \ldots, \beta_t)$ satisfying above conditions. We shall only consider the case when
$s,t>0$ because the other three cases are trivial.

Set $\bar{\alpha}_i=\alpha_i-m+t$ and $\bar{\beta}_j=\beta_j-k-s-j+1$.
It is easy to chcek that $\bar{\alpha}=(\bar{\alpha}_1,
\bar{\alpha}_2, \ldots, \bar{\alpha}_s)$ and
$\bar{\beta}=(\bar{\beta}_1, \bar{\beta}_2, \ldots, \bar{\beta}_t)$
form two partitions with $\bar{\alpha}_s \geq t$ and
$\bar{\beta}_t \geq 0$. So we may insert $\bar{\beta}$ into $\bar{\alpha}$
to create a pair of partitions $(\mu,\nu)$  via the following
procedure.

\begin{itemize}

\item[$\bullet$ ] If $\bar{\beta}_1 \geq \bar{\alpha}_1$, we insert
$\bar{\beta}_1$ into $\bar{\alpha}$ to form a new partition
$\delta=(\delta_1, \delta_2, \ldots, \delta_{s+1})=(\bar{\beta}_1,
\bar{\alpha}_1-1, \bar{\alpha}_2-1, \ldots, \bar{\alpha}_s-1)$.
Moreover, we set $\nu_1=k+s+t-1$ to  record the insertion position.

\item[$\bullet$ ] Otherwise, we assume that  $j_1$ is the largest
integer such that $\bar{\alpha}_{j_1} > \bar{\beta}_1 $. Then we
insert $\bar{\beta}_1$ into $\bar{\alpha}$ to form a new partition
$\delta=(\delta_1, \delta_2, \ldots, \delta_{s+1})=(\bar{\alpha}_1,
\ldots, \bar{\alpha}_{j_1}, \bar{\beta}_1, \bar{\alpha}_{j_1+1}-1,
\ldots, \bar{\alpha}_{s}-1)$. In this case, we use
$\nu_1=k+s+t-j_1-1$ to record the insertion position of $\bar{\beta}_1$.
Obviously, $k+t-1 \leq \nu_1 \leq k+s+t-2$.

\end{itemize}

Conversely, given a partition $\delta=(\delta_1, \delta_2, \ldots,
\delta_{s+1})$ and a number $\nu_1$ with $k+t-1 \leq \nu_1 \leq
k+s+t-1$, we may extract the part $\bar{\beta}_1$ from $\delta$. It
is clear that the above procedure is reversible.

Similarly, we can insert $\bar{\beta}_2$ into the
partition $\delta=(\delta_1, \delta_2, \ldots, \delta_{s+1})$.
Applying the insertion algorithm repeatedly to $\bar{\beta}_2,
\ldots, \bar{\beta}_t$, we come to a partition $\bar{\mu}$ with
$s+t$ parts and the desired partition $\nu$, where $\nu=(\nu_1,
\ldots, \nu_t)$ with $\nu_i=k+s+t-j_i-1$ for each $1 \leq i \leq t$.
Furthermore, one sees that $k+t-i \leq \nu_i \leq k+s+t-i$. On the
other hand, we get the desired partition $\mu$ by setting
$\mu=\{\bar{\mu}_1+m,\ldots,\bar{\mu}_s+m,\bar{\mu}_{s+1},\ldots
\bar{\mu}_{s+t}\}$. This completes the proof.  \qed

 In the above correspondence, the partition $\mu$ is also called  the insertion
partition and $\nu$ is called the record partition.
 As an example, let $k=3,m=2$,$s=4, t=3$ and $\alpha=(8, 7,5,3)$,
$\beta=(12,11,9)$. Then we have $\bar{\alpha}=(9,8,6,4)$ and
$\bar{\beta}=(5,3,0)$, and
$$\bar{\mu}=(9,8,6,5,3,2,0), \quad \nu=(6,5,3),\quad \mu=(11,10,8,7,3,2,0).$$
The above correspondence is illustrated as follows
\begin{center}
\begin{tabular}{ccccccccc}
9 & 8 & 6 &  &  &4 \\
\hline 9 & 8 & 6 & 5 & 3 & 2& 0\\
  & & & 6& 5& &3 \\
  \hline
 & & & 5&3 &  &0\\
\end{tabular}
\end{center}

It is worth mentioning that the conditions $\alpha_s>m$ and
$\beta_t > k+s+t-1$ can be recast in terms of the single  statement
that $\mu_{s+t}>0$. This observation will be useful in the proof of
Theorem \ref{66}.

We now turn our attention to  the minor difference between the Algorithm
Z and the above variation. Given two partitions $\alpha$ and $\beta$, we may
apply the
Algorithm Z to a pair of partitions $(\bar{\alpha},\bar{\beta})$,
where
$\bar{\alpha}=(\alpha_1-s-m+1,\alpha_2-s-m+2,\ldots,\alpha_s-m)$ and
$\bar{\beta}=(\beta_1-s-t-k+1,\ldots,\beta_t-s-t-k+1)$. It can be
seen that the record partition of $(\alpha,\beta)$ and the record
partition of $(\bar{\alpha},\bar{\beta})$ differ only by  a
staircase partition $(k+t-1,k+t-2,\ldots,k)$. For the above example,
one has $\bar{\alpha}=(3,3,2,1)$, $\bar{\beta}=(3,2,0)$. Inserting
$\bar{\beta}$ into $\bar{\alpha}$ via the Algorithm Z gives
$\bar{\mu}=(3,3,2,2,1,1,0)$, $\bar{\nu}=(1,1,0)$ as illustrated below
\begin{center}
\begin{tabular}{ccccccccc}
3 & 3 & 2 &  &  &1 \\
\hline 3 & 3 & 2 & 2 & 1 & 1& 0\\
  & & & 1& 1& &0 \\
  \hline
 & & & 3&2 &  &0\\
\end{tabular}
\end{center}

It is not hard to see  that our combinatorial proof of Ramanujan's
formula can be restated in terms of the original Algorithm Z. Nevertheless,
the variation seems to be more convenient for the sake of presentation.

\begin{coro}
There is a bijection between the set of pairs of partitions
$(\alpha, \beta)$ and  the set of triples of partitions $(n; \mu,
\nu, \gamma)$ satisfying the conditions
\begin{itemize}
\item[${\bullet}$] $\alpha$ has distinct nonnegative parts and $\beta$ has
    nonnegative
parts;
\item[${\bullet}$] $\mu$ has $n$ distinct nonnegative parts,
$\nu$ has either distinct nonnegative parts with each part $\leq
n-1$ (corresponding to $\beta_1\geq l(\alpha)$) or is an empty
partition (corresponding to $\beta_1<l(\alpha)$), and $\gamma$ has
nonnegative parts with each part $\leq n-1$;
\item[$\bullet$ ] $|\alpha|+|\beta|=|\mu|+|\nu|+|\gamma|$.
\end{itemize}
\end{coro}

\pf  Assume that  $n$ is the largest number satisfying $\beta_{n-l(\alpha)}\ge n-1$. If such an $n$ exists, then set $\gamma=(\beta_{n-l(\alpha)+1},\ldots,\beta_{l(\beta)})$, which is a partition
with each part $\le n-1$.  Denote
by $\bar{\beta}$ the partition
$(\beta_1,\ldots,\beta_{n-l(\alpha)})$and apply the bijection
$\varphi$ in Theorem \ref{7} to
$(\alpha, \bar{\beta})$ for $m=0$ and $k=0$, we
get a partition $\mu$  having $n$ distinct nonnegative parts and a
partition $\nu$ having distinct parts with $0\leq \nu_i \leq n-1$.
If there does not exist such an $n$, namely $\beta_1 \leq
l(\alpha)-1$,  then we set $n=l(\alpha)$, $\mu =\alpha$,
$\gamma=\beta$ and set $\nu$ to be the empty partition.
This completes the proof. \qed

The above corollary can be regarded as a combinatorial
interpretation of the following identity \cite[Exercise 1.6
(ii)]{GasperRahman}:
\begin{equation}\label{11}
\frac{(-a;q)_\infty}{(b;q)_\infty}=\sum_{n=0}^\infty
\frac{P_n(a,-b)q^{\binom{n}{2}}}{(q;q)_n(b;q)_n}.
\end{equation}

\section{The Combinatorial Proof} \label{S:P*}

In this section, we  aim to give a combinatorial proof of Ramanujan's ${}_1\psi_1$ summation formula
\eqref{32}. When $N \geq 0$, the coefficient of $z^N$  on the
left-hand side equals the generating function for the quintuples $(n;
\alpha, \beta, \gamma, \lambda, \mu)$ subject to the following conditions:
  \begin{itemize}
\item[${\bullet}$] $\alpha$ has distinct and positive parts,
\item[${\bullet}$] $\beta$ has positive parts,
\item[${\bullet}$] $\gamma$ has distinct nonnegative parts,
\item[${\bullet}$] $\lambda$ has distinct nonnegative parts with each part $\leq n-1$,
\item[${\bullet}$] $\mu$ has nonnegative parts with each part $\leq n-1$,
\end{itemize}
where the exponents of $a$ and $b$ are used to keep track  of
$l(\lambda)-l(\alpha)-l(\gamma)$ and $l(\gamma)+l(\mu)$
respectively, and $N$ records $n-l(\gamma)$.
 The coefficient of $z^N$ on the right-hand side is the generating function
for the quintuples $(A, B, C, D, E)$ of partitions with the following restrictions:
\begin{itemize}
\item[${\bullet}$] Both $A$ and $C$ have distinct nonnegative parts,
\item[${\bullet}$] Both $B$ and $D$ have nonnegative parts,
\item[${\bullet}$] $E$ has distinct and positive parts,
\end{itemize}
where the exponents of $a$ and $b$ are used to keep track of
$l(C)-l(A)-l(E)$ and $l(A)+l(B)$ respectively, and  $N$  records the number $l(C)+l(D)-l(E)$.
Let $\mathscr{A}$ and  $\mathscr{B}$ be the sets of the quintuples
$(n;\alpha, \beta, \gamma, \lambda, \mu)$ and  $(A, B, C, D, E)$,
as defined above.
\begin{theo}\label{66}
There is a bijection between $\mathscr{A}$ and $\mathscr{B}$.
\end{theo}

\pf Given a quintuple $(n; \alpha, \beta, \gamma, \lambda, \mu)$
with $N=n-l(\gamma)$. As an example, for $N=4$, $n=9$, let
\[ \alpha=(10,9,5,3,2),\quad  \beta= (13,11,10,9,9,5,4,4,2), \quad \gamma=(9, 6,
4, 2, 1),\] \[ \lambda=(7, 6, 5, 3, 1), \quad \mu=(5, 4, 4, 1).\] We shall use
this example to illustrate the operations at every step.

\noindent
Step 1.
Find the largest number $p$ such that $\lambda_p\ge N$. Then
$\bar{\lambda}=(\lambda_1,\ldots,\lambda_p)$ is a partition into
distinct parts with $N \leq \lambda_i \le n-1$ and $1\leq i \leq p$.
Let $F=(\lambda_{p+1},\ldots,\lambda_{l(\lambda)})$, where $p$ is the
largest number such that $\lambda_{p+1} \leq N-1$. Clearly, $l(\bar{\lambda})\leq n-N=l(\gamma)$.

Applying the bijection $\varphi^{-1}$ in Theorem
\ref{7} to the pair ($\gamma$, $\bar{\lambda}$) with $m=0$ and
$k= N$, we obtain a partition $A$ with distinct
nonnegative parts, and a partition $\bar{B}$ with every  part $\ge
n-1$. Now we can put the parts of  $\bar{B}$ and $\mu$ together to form
the desired partition $B$. Note that if such an integer $p$ does not
exist, that is, $\lambda_1\leq N-1$, then we have $F = \lambda$, $A =
\gamma$ and $B = \mu$.

 For  the above example, we have
 \[ \bar{\lambda}=(7,6,5),
\quad F=(3,1), \quad  A=(6, 1), \quad B=(12, 11, 10, 5, 4, 4, 1).\]

\noindent Step 2.
Find the largest number $l$ such that $\beta_{l-l(\alpha)}\ge N+l$,
and set $\bar{\beta}=(\beta_1,\ldots,\beta_{l-l(\alpha)})$. Add
enough zero parts if necessary to the conjugate of the partition
$(\beta_{l-l(\alpha)+1},\ldots,\beta_{l(\beta)})$ to obtain a
partition $\bar{D}$ with   $N+l$ parts.

Now we can apply the mapping  $\varphi$ to ($\alpha$,
$\bar{\beta}$) with $m=0$ and $k =N$ to generate
a partition $E$ with $l$ distinct positive parts, since
$\alpha_{l(\alpha)}>0$ and $\beta_{l-l(\alpha)}\geq N+l$. Meanwhile, we also
obtain a
partition $\bar{F}$ into distinct parts with each part $\ge N$
$\le N+l-1$. So we can put  $F$ and $\bar{F}$ together to create a partition
$\bar{C}$ into distinct nonnegative parts with each part $\leq
N+l-1$. Note that if such an integer $l$ does not exist, that is,
$\beta_1 \leq N+l(\alpha)$, we may set $E =\alpha$ and $\bar{C}=F$. In this
case, $\bar{D}$ is a partition with $N+l(\alpha)$ parts, which can be
obtained from the conjugate of $\beta$ with some zero parts added  if
needed.

 For the above example, we have $$\bar{\beta}=(13,11), \quad
\bar{D}=(7,7,6,6,4,3,3,3,3,1, 0),$$
$$
E=(12,11,7,5,4,3,1), \quad  \bar{C}=(6,4,3,1).$$

Applying the bijection $\phi^{-1}$ in Corollary \ref{4}
to ($\bar{D}$,$\bar{C}$), we obtain  the  partition $C$
 into distinct nonnegative parts and the partition  $D$
into nonnegative parts. For the above example,
we have
\[ C=(10, 7, 6, 2), \quad  D=(7, 7, 6, 6, 3, 3,
0).\]

Whence we have constructed a quintuple $(A, B, C, D, E)$ for which
\[ |A|+|B|+|C|+|D|+|E|=|\alpha|+|\beta|+|\gamma|+|\lambda|+|\mu|.\]
Notice that the exponents of $a$ and $b$ remain unchanged during
the above procedure. Since each step is reversible, we have
established a bijection between $\mathscr{A}$ and $\mathscr{B}$.
This completes the proof. \qed

When $N=-m < 0$,  by multiplying  both sides of \eqref{32} by
$\frac{(b;q)_{-m}}{(-a;q)_{-m}}$, we get
\begin{multline}
\frac{(-q/a;q)_{\infty}(-b/az;q)_{\infty}}{(q;q)_{\infty}}\sum_{l=0}^{\infty}\frac{(-aq^{-m};q)_l}{(bq^{-m};q)_l}z^{l-m}\\
 = \frac{(-b/a;q)_{\infty}(-az;q)_\infty(-q/az;q)_\infty (-aq^{-m};q)_m}{ (bq^{-m};
q)_\infty(z;q)_\infty}.
 \end{multline}
Substituting $b$ by $bq^m$ and using  Euler's identity,
$$(-bq^m/az;q)_{\infty}=\sum_{n=0}^{\infty}\frac{(b/az)^nq^{mn+{n \choose 2}}}{(q;q)_n},$$
the coefficients of $z^N$ on both sides can be written as
\begin{multline} \label{eq:36}
\frac{a^{-m}q^{m+1\choose 2}(-q^{m+1}/a;
q)_{\infty}}{(q;
q)_{\infty}}\sum_{l=0}^{\infty}\frac{P_l(bq^m/a,
-b)q^{l\choose
2}}{(b;q)_l(q;q)_l}\\
 = \big[z^{-m} \big] \frac{(-bq^m/a;q)_{\infty}(-az;q)_\infty(-q/az;q)_\infty}{ (b;
q)_\infty(z;q)_\infty},
 \end{multline}
where $[x^n]F(x)$ denotes the coefficient of $x^n$ in $F(x)$.

 Each term on the left-hand side of \eqref{eq:36} can be interpreted as the
generating function for the quintuples $(l;\alpha, \beta, \gamma,
\lambda,\mu)$ defined as follows
\begin{itemize}
\item[${\bullet}$] $\alpha$ has distinct and positive
parts with $\alpha_{l(\alpha)-m+1}=m$,
\item[${\bullet}$] $\beta$ has positive parts,
\item[${\bullet}$] $\lambda$ has  distinct nonnegative parts with each part $\leq l-1$,
\item[${\bullet}$] $\gamma$ has distinct $l$
positive parts. Let $s=l-l(\lambda)$. Then $ \gamma_{s}-\gamma_{s+1}\ge
m+1$ if $0 <l(\lambda) <l$ and
$\gamma_s \geq m$ if $l(\lambda)=0$.
\item[${\bullet}$] $\mu$ has nonnegative parts with each part $\leq l-1$,
\end{itemize}
where the exponent of $a$ records $l(\lambda)-l(\alpha)-l(\gamma)$
and the exponent of $b$ keeps track of $l(\gamma)+l(\mu)$.

Clearly,
the right-hand side of \eqref{eq:36} is the generating function for
 the quintuples $(A, B, C, D, E)$ defined as follows
\begin{itemize}
\item[${\bullet}$] $A$ has distinct parts with each part $\geq m$,
\item[${\bullet}$] Both $B$ and $D$ have nonnegative parts,
\item[${\bullet}$] $C$ have distinct nonnegative parts,
\item[${\bullet}$] $E$ has distinct and positive parts,
\end{itemize}
where the exponent of $a$ records
 $l(C)-l(A)-l(E)$, the exponent
   of $b$ keeps track of  $l(A)+l(B)$
 and $l(C)+l(D)-l(E)=-m$.

Let $\mathscr{C}$ and $\mathscr{D}$ be the sets of  quintuples
$(l;\alpha, \beta, \gamma, \lambda, \mu)$ and $(A, B, C, D, E)$
as given before.

\begin{theo}
There is a bijection between $\mathscr{C}$ and $\mathscr{D}$.
\end{theo}

\pf  Let  $(A, B, C, D, E)$ be a quintuple with $N=l(C)+l(D)-l(E)$.
As an example, for $N=-2$, assume that
 \[ A=(12,11,7,5,4),\quad
B=(15,13,12,11,11,7,6,6,4,2,1,1),\quad  C=(7,6,4,1,0),\]
\[ D=(9,8,5,5,4,1),\quad  E=(22,19,18,17,15,12,11,10,8,7,6,3,1).\]
We shall use this example to illustration the operation at each step.

\noindent
Step 1. If $B_1\leq l(A)-1$,
we  set $l=l(A)$, $\mu=B$, $\gamma=A$ and set $\lambda$ to be the
empty partition. Otherwise, we find the largest number $l$ such that
$B_{l-l(A)}\ge l-1$, and set $\bar{B}=(B_1,\ldots,B_{l-l(A)})$. Now
$(B_{l-l(A)+1},\ldots,B_{l(B)})$ is  the desired partition $\mu$.
Apply  $\varphi$ in Theorem \ref{7} to ($A$, $\bar{B}$)
with $k=0$. In the case $s=l(A)>0$, we
get a partition $\gamma$ into $l$ distinct nonnegative parts with
$\gamma_{l(A)}-\gamma_{l(A)+1}\ge m+1$ and a partition $\lambda$
into distinct nonnegative parts with each part $\le l-1$. When
$s=0$, we get a partition $\gamma=(B_1,B_2-1,\ldots,B_l-l+1)$ and a
partition $\lambda=(l-1,l-2,\ldots,0)$.

 For the above  example, we have
$$ \bar{B}=(15,13,12,11,11), \mu=(7, 6,
6, 4, 2,1,1),$$
$$\gamma=(17,16,12,11,9,6,5,4,3,2), \lambda=(7,5,3,1,0).$$

\noindent
Step 2.
Applying the bijection $\phi$ in Corollary \ref{4} to $(C, D)$
yields a partition $\bar{C}$ into distinct nonnegative parts with
 each part
$\le n-1$ and a partition $\bar{D}$
 into $n$ nonnegative parts. Evidently, we have $l(C)+l(D)=n$, and
 hence $l(E)=n+m$.

 For the above example, we find
$\bar{C}=(5,4,3,1,0)$,
$\bar{D}=(9,8,5,5,4,2,2,1,1,0,0)$.

\noindent
Step 3.
After removing a staircase partition $(n+m,n+m-1,\ldots,1)$ from
$E$, we are left with  a partition with $l(E)$ nonnegative parts,
whose conjugate is denoted by $\bar{E}$. Add $m+1$ to each part of
$\bar{C}$ to obtain a partition $\tilde{C}$. Now we may construct a
partition  $\tilde{D}$ by adding a staircase $(n-1,n-2,\ldots,0)$ to
$\bar{D}$, then adding $m+1$ to the first $l(D)$ parts.

 Applying $\varphi^{-1}$ in Theorem \ref{7} to
($\tilde{D}$, $\tilde{C}$) with $k$ replaced by $m+1$ and
 $m$ replaced by $m+1$ yields
 a partition $\bar{\alpha}$ into $l(D)$ distinct parts with
each part $\ge m+1$ and a partition $\bar{\beta}$ into $l(C)$ parts
with each part $\ge n+m=l(E)$. Combining $\bar{\alpha}$ with a
staircase partition $(m,m-1,\ldots,1)$ gives the  partition
$\alpha$, and combining $\bar{\beta}$ with $\bar{E}$ gives the
required partition $\beta$.

For the above example, we get
$$\bar{E}=(12,11,11,8,5,5,4,1,1), \quad \tilde{C}=(8,7,6,4,3), \quad
\tilde{D}=(22,20,16,15,13,10,6,4,3,1,0),$$
$$\beta=(16,16,15,13,13,12,11,11,8,5,5,4,1,1),\quad
\alpha=(17,15,11,10,8,4,2,1).$$

Thus we have constructed a quintuple $(l;\alpha, \beta, \gamma, \lambda,
 \mu)$ such that
 $$|\alpha|
+|\beta|+|\gamma|+|\lambda|+|\mu|=|A|+|B|+|C|+|D|+|E|.$$
Moreover, the exponents of $a$ and $b$ are preserved
at every step. It should be mentioned that  each step of the above procedure is
reversible. This completes the proof. \qed

{\bf Acknowledgments.} This work was supported by the 973 Project,
the PCSIRT Project of the Ministry of Education,  and the National
Science Foundation of China.

\end{document}